\documentclass[12pt,fleqn]{article}

\usepackage{amsmath}
\usepackage{amsthm}
\usepackage{amssymb}
\usepackage{amsfonts}
\usepackage{graphicx}
\usepackage{color}
\newtheorem{theorem}{Theorem}[section]

\newtheorem{corollary}[theorem]{Corollary}

\theoremstyle{definition}
\newtheorem{definition}[theorem]{Definition}

\newtheoremstyle{named}{}{}{\itshape}{}{\bfseries}{.}{.5em}{\thmnote{#3's }#1} \theoremstyle{named} 

\theoremstyle{remark}
\newtheorem{remark}[theorem]{Remark}

\numberwithin{equation}{section}

\numberwithin{equation}{section}

\title{Existence of a lower bound for the distance between point masses of relative equilibria for generalised quasi-homogeneous $n$-body problems and the curved $n$-body problem}
\author{Pieter Tibboel$^\ast$\\
Date of resubmission: The 8th of February, 2015\date{ }
}

\begin{document}
\maketitle
\begin{abstract}
  We prove that if for relative equilibrium solutions of a generalisation of quasi-homogeneous $n$-body problems the masses and rotation are given, then the minimum distance between the point masses of such a relative equilibrium has a universal lower bound that is not equal to zero. We furthermore prove that the set of such relative equilibria is compact and prove related results for $n$-body problems in spaces of constant Gaussian curvature. 
\end{abstract}

\begin{description}

\item \hspace*{3.8mm}$\ast$ Department of Mathematical Sciences, Xi'an Jiaotong-Liverpool University, Suzhou, China. \\
Email: \texttt{Pieter.Tibboel@xjtlu.edu.cn}

\end{description}

\newpage
\section{Introduction}
By $n$-body problems we mean problems where we are tasked with deducing the dynamics of $n$ point masses. The study of such problems has obvious applications to fields such as celestial mechanics, chemistry, atomic physics and crystallography (see for example \cite{AbrahamMarsden}, \cite{CDL}, \cite{CLP}, \cite{CRS}, \cite{D-3}, \cite{D-1}, \cite{D0}, \cite{DDLMMPS}, \cite{DMS}, \cite{DPS1}, \cite{DPS2}, \cite{DPS3}, \cite{DPS4}, \cite{PSSY}, \cite{PV} and the references therein). The $n$-body problems discussed in this paper are the $n$-body problem in spaces of constant Gaussian curvature and a generalisation of quasi-homogeneous $n$-body problems, which we will call generalised quasi-homogeneous $n$-body problems for short:
\begin{definition}\label{Definition n-body problem in spaces of constant curvature}
  Let   \begin{align*}
    \mathbb{M}_{\sigma}^{k}=\{(x_{1},....,x_{k+1})\in\mathbb{R}^{k+1}|x_{1}^{2}+...+x_{k}^{2}+\sigma x_{k+1}^{2}=\sigma\},
  \end{align*}
  where $\sigma$ equals either $+1$, or $-1$
  and for $x$, $y\in\mathbb{M}_{\sigma}^{k}$ define the inner product
  \begin{align*}
    x\odot y=x_{1}y_{1}+...+x_{k}y_{k}+\sigma x_{k+1}y_{k+1}.
  \end{align*}
  By the $n$-body problem in spaces of constant Gaussian curvature, henceforth referred to as \textit{the $n$-body problem in spaces of constant curvature}, or \textit{curved $n$-body problem} (see \cite{DPS1}, \cite{DPS2} and \cite{DPS3}), we mean the problem of finding the dynamics of $n$ point particles with respective masses $m_{1}$,..., $m_{n}$ and coordinates $q_{1}$,..., $q_{n}\in\mathbb{M}_{\sigma}^{k}$, $k\geq 2$, as described by the system of differential equations
\begin{align}\label{EquationsOfMotion Curved}
   \ddot{q}_{i}=\sum\limits_{j=1,\textrm{ }j\neq i}^{n}\frac{m_{j}(q_{j}-\sigma(q_{i}\odot q_{j})q_{i})}{(\sigma-\sigma(q_{i}\odot q_{j})^{2})^{\frac{3}{2}}}-\sigma(\dot{q}_{i}\odot\dot{q}_{i})q_{i},\textrm{ }i\in\{1,...,n\}.
\end{align}
\end{definition}
\begin{definition}\label{Definition generalised n-body problem}
  By \textit{generalised quasi-homogeneous $n$-body problems} we mean problems where we have to find the dynamics of $n$ point particles $q_{1}$,..., $q_{n}\in\mathbb{R}^{k+1}$, $k\geq 1$ with respective masses $m_{1}$,...,$m_{n}$ as described by the system of differential equations \begin{align}\label{Equations of motion}
  \ddot{q}_{i}=\sum\limits_{j=1,\textrm{ }j\neq i}^{n}m_{j}(q_{j}-q_{i})f(\|q_{j}-q_{i}\|),\textrm{ }i\in\{1,...,n\},
\end{align}
where $f:\mathbb{R}_{>0}\rightarrow\mathbb{R}$ can be any continuous function with the property that $f(x)$ and $xf'(x)$ are bounded and differentiable for $x$ away from $0$ and $\lim\limits_{x\downarrow 0}f(x)=\pm\infty$. If $f(x)=x^{-\frac{3}{2}}$ and $k=2$, then we speak of the \textit{classical $n$-body problem}. If $f(x)=ax^{-\alpha}+bx^{-\beta}$, $\alpha$, $\beta\in\mathbb{R}_{>0}$, $a$, $b\in\mathbb{R}$, then we speak of a \textit{quasi-homogeneous $n$-body problem}. The reason that $xf'(x)$ needs to be bounded for $x$ away from zero is a technical one: In the proof of Theorem~\ref{limits distinct}, $\|Q_{2r}-Q_{jr}\|B_{jr}$ needs to be bounded. $B_{jr}$ depends on $f'$ and $\|Q_{2r}-Q_{jr}\|$ may be unbounded. $xf'(x)$ being bounded for $x$ away from zero is sufficient to ensure that $\|Q_{2r}-Q_{jr}\|B_{jr}$ is bounded and does not impose very strong restrictions on the generality of our problem.
\end{definition}
Research into $n$-body problems for spaces of constant Gaussian curvature goes back as far as the 1830s, when Bolyai and Lobachevsky (see \cite{BB} and \cite{Lo} respectively) independently proposed a curved 2-body problem in hyperbolic space $\mathbb{H}^{3}$. Since then, $n$-body problems in spaces of constant Gaussian curvature have been investigated by mathematicians such as Dirichlet, Schering (see \cite{S1}, \cite{S2}), Killing (see \cite{K1}, \cite{K2}, \cite{K3}), Liebmann (see \cite{L1}, \cite{L2}, \cite{L3}) and Kozlov and Harin (see \cite{KH}). However, the succesful study of $n$-body problems in spaces of constant Gaussian curvature for the case that $n\geq 2$ began with \cite{DPS1}, \cite{DPS2}, \cite{DPS3}  by Diacu, P\'erez-Chavela and Santoprete. After this breakthrough, further results for the $n\geq 2$ case were then obtained in \cite{CRS}, \cite{D1}, \cite{D2}, \cite{D3}, \cite{D4}, \cite{D5}, \cite{DK}, \cite{DP}, \cite{DPo}, \cite{DT}, \cite{T}, \cite{T2} and \cite{T3}. For a more detailed historical overview, please see \cite{D2}, \cite{D3}, \cite{D4}, \cite{D6}, \cite{DK}, or \cite{DPS1}.

Quasi-homogeneous $n$-body problems for general values of $n$ started with \cite{D0} by Diacu and can be applied to many fields ranging from celestial mechanics to atomic physics to chemistry to crystallography (see \cite{DPS4}). For examples see  \cite{CDL}, \cite{CLP}, \cite{D-3}, \cite{D-1}, \cite{D0}, \cite{DDLMMPS}, \cite{DMS}, \cite{Kuzmina}, \cite{Llibre}, \cite{Moulton}, \cite{Palmore}, \cite{PSSY}, \cite{PV}, \cite{Roberts}, \cite{Wintner} and the references therein.

Solutions to $n$-body problems that constitute point configurations that retain their size and shape over time are called \textit{relative equilibria}. However, following \cite{D3}, we will use a slightly more general definition: In Euclidean space, any configuration $q_{1}(t)$,...,$q_{n}(t)$ that retains its size and shape over time can be written as $q_{1}(t)=T(t)Q_{1}$,...,$q_{n}(t)=T(t)Q_{n}$, where $Q_{1}$,...,$Q_{n}\in\mathbb{R}^{k}$ are constant vectors and $T(t)$ is a a block matrix of rotation matrices. Note that for any rotation matrix $T(t)$ and $x$, $y\in\mathbb{R}^{k}$, we have that $\langle T(t)x, T(t)y\rangle=\langle x,y\rangle$, as rotations preserve the distance between points and angles between lines. Consequently, all matrices $T(t)$ for which $\langle T(t)x, T(t)y\rangle=\langle x,y\rangle$ are norm preserving (take $x=y$) and thus preserve angles (as for any vectors $x$, $y$, the norms of $x$, $y$ and $x-y$ are preserved and thus the angles of the triangle with sides of length $\|x\|$, $\|y\|$ and $\|x-y\|$ are preserved). On $\mathbb{R}^{k}$, $T(t)$ being any block matrix of rotation matrices is equivalent with $\langle T(t)x, T(t)y\rangle=\langle x,y\rangle$. On $\mathbb{M}_{\sigma}^{k}$, if we look for $T(t)$ for which $x\odot_{k}y=(T(t)x)\odot_{k}(T(t)y)$, $T(t)$ that depend on hyperbolic and parabolic functions exist as well (see \cite{D3}). As we do not explicitly use the many different possible expressions for $T(t)$, it makes sense to therefore define relative equilibria as follows:
\begin{definition}\label{Definition Relative Equilibrium}
  Consider any solution to (\ref{EquationsOfMotion Curved}), or (\ref{Equations of motion}) for which the $q_{1}$,..., $q_{n}$ can be written as  $q_{1}(t)=T(t)Q_{1}$,..., $q_{n}(t)=T(t)Q_{n}$, where $Q_{1}$,..., $Q_{n}\in\mathbb{R}^{k+1}$ are constant and $T(t)$ is a time dependent, invertible $k\times k$ matrix for which there exist constants $c_{1}$, $c_{2}\in\mathbb{R}_{>0}$ such that for all $x\in\mathbb{R}^{k}$ $c_{1}\|x\|\leq\|T(t)^{-1}\ddot{T}(t)x\|\leq c_{2}\|x\|$. If for all $x$, $y\in\mathbb{R}^{k}$ we have that $\langle T(t)x,T(t)y\rangle=\langle x,y\rangle$, we call solutions of this type that solve (\ref{Equations of motion}) a \textit{relative equilibrium} of (\ref{Equations of motion}). If for all $x$, $y\in \mathbb{M}_{\sigma}^{k}$ we have that  $(T(t)x)\odot(T(t)y)=x\odot y$ and $\|\dot{T}(t)x\|$ is bounded if $x$ is bounded, then we call solutions of this type that solve (\ref{EquationsOfMotion Curved}) a \textit{relative equilibrium} of (\ref{EquationsOfMotion Curved}). For generalised quasi-homogeneous $n$-body problems, we call the set of all such configurations that are equivalent under rotation and scalar multiplication a \textit{class of relative equilibria}. It should be remarked that the property that $c_{1}\|x\|\leq\|T(t)^{-1}\ddot{T}(t)x\|\leq c_{2}\|x\|$ does not follow directly from the fact that $x\odot_{k}y=(T(t)x)\odot_{k}(T(t)y)$. Checking the different possible cases in \cite{D3} shows easily that the given properties of $T(t)$ are true and saves us many pages of writing out block matrices.
\end{definition}
Relative equilibria can tell a great deal about the physical space for which their respective $n$-body problems have been defined: A prime example of how much information can be deduced by studying relative equilibria, comes from celestial mechanics: It was proven in \cite{DPS1} and \cite{DPS2} that for the $n$-body problem in spaces of constant curvature (i.e. spheres or hyperboloids) relative equilibria that are shaped as equilateral triangles have to have equal masses. This means that our solar system, with the Sun, Jupiter and the Trojan asteroids forming approximately an equilateral triangle and relative equilibrium, is likely flat within the respective area, i.e. has zero Gaussian curvature. For further information on the relevance of relative equilibria, see for example \cite{DPS0} and \cite{Saari}. Of particular importance is the link with the sixth Smale problem (see \cite{Smale}), which states that for the classical case, if the equilibria are induced by a plane rotation, the number of classes of relative equilibria is finite, if the masses $m_{1}$,...,$m_{n}$ are given. This problem is still open for $n>5$ and was solved for $n=3$ by A. Wintner (see \cite{Wintner}), $n=4$ by M. Hampton and R. Moeckel (see \cite{HamptonMoeckel}) and for $n=5$ by A. Albouy and V. Kaloshin, assuming that the 5-tuple of positive masses belongs to a given codimension 2 subvariety of the mass space (see \cite{AlbouyKaloshin}).
As a potential step towards a proof of Smale's problem, M. Shub showed in \cite{Shub} that the set of all classes of relative equilibria, provided they have the same set of masses, is compact. Additionally, Shub proved in the same paper that if the rotation inducing the equilibria is given as well, there exists a universal nonzero, minimal distance that the point masses lie apart from each other.

In this paper, as a logical next step after Shub's work in \cite{Shub} and to gain further understanding of the geometry of relative equilibria, we prove Shub's results when using (\ref{Equations of motion}) instead of the classical $n$-body problem and related results for $n$-body problems in spaces of constant curvature. Specifically, we prove that
\begin{theorem}\label{limits distinct}
  Consider the set $R_{T,m_{1},...,m_{n}}$ of all relative equilibria of (\ref{Equations of motion}) with rotation matrix $T(t)$ and masses $m_{1}$,..., $m_{n}$. Then there exists a constant $c\in\mathbb{R}_{>0}$ such that for all relative equilibria $\{T(t)Q_{i}\}_{i=1}^{n}$ in the set $R_{T,m_{1},...,m_{n}}$, we have that $\|Q_{i}-Q_{j}\|>c$ for all $i$, $j\in\{1,...,n\}$, $i\neq j$.
\end{theorem}
and consequently that
\begin{corollary}\label{compact}
  Consider the set $R_{T,m_{1},...,m_{n}}$ of all relative equilibria (\ref{Equations of motion}) with rotation matrix  $T(t)$ and masses $m_{1}$,..., $m_{n}$. If $\lim\limits_{x\downarrow 0}xf(x)=\pm\infty$ and $xf(x)$ is bounded for $x$ away from $x=0$, then there exists a $C\in\mathbb{R}_{>0}$ such that for all relative equilibria $\{T(t)Q_{i}\}_{i=1}^{n}$ in the set $R_{T,m_{1},...,m_{n}}$, we have that $\|Q_{i}\|<C$ for all $i\in\{1,...,n\}$.
\end{corollary}
%
For the $n$-body problem in spaces of constant curvature, proving that relative equilibria form a compact set is pointless for the case that $\sigma=1$ (i.e. the problem is defined on the unit sphere) but it does make sense to investigate whether there exists a universal lower bound for the distance between the point masses. To that extent, we will prove that
\begin{theorem}\label{Main Theorem Curved Case 1}
  Let $\sigma=1$, $\epsilon>0$ and let $R_{\epsilon, T,m_{1},...,m_{n}}$ be the set of all relative equilibria $T(t)Q_{1}$,...,$T(t)Q_{n}$ of (\ref{EquationsOfMotion Curved}) for which $\langle Q_{i},Q_{j}\rangle>-1+\epsilon$, $i$, $j\in\{1,...,n\}$, $i\neq j$. Then for all $\epsilon>0$ there exists a constant $c_{\epsilon}>0$ such that for any relative equilibrium solution in $R_{\epsilon, T,m_{1},...,m_{n}}$ of (\ref{EquationsOfMotion Curved}), $\|q_{i}-q_{j}\|_{k}>c_{\epsilon}$ for all $i$, $j\in\{1,...,n\}$, $i\neq j$ if the masses $m_{1}$,..., $m_{n}$ and rotation $T(t)$ are given.
\end{theorem}
For the $n$-body problem on a hyperbola, relative equilibria do not form a compact set with respect to the Euclidean norm. For the negative curvature case, we will prove that
\begin{theorem}\label{Main Theorem Curved Case 2}
  Let $\sigma=-1$ and let $R_{T,m_{1},...,m_{n}}$ be the set of all relative equilibria $T(t)Q_{1}$,...,$T(t)Q_{n}$ of (\ref{EquationsOfMotion Curved}). Then for any bounded subset $W$ of $R_{T,m_{1},...,m_{n}}$ there exists a constant $C_{W}>0$ such that for any relative equilibrium solution in $W$ of (\ref{EquationsOfMotion Curved}), $\|q_{i}-q_{j}\|_{k}>C_{W}$ for all $i$, $j\in\{1,...,n\}$, $i\neq j$ if the masses $m_{1}$,..., $m_{n}$ and rotation $T(t)$ are given.
\end{theorem}
\begin{remark}
  Theorem~\ref{Main Theorem Curved Case 1} and Theorem~\ref{Main Theorem Curved Case 2} were proven for a very specific subclass of relative equilibria in \cite{T3}. In this paper, it should be noted that Theorem~\ref{Main Theorem Curved Case 1} and Theorem~\ref{Main Theorem Curved Case 2} hold for all types of relative equilibria (positive elliptic relative equilibria, positive elliptic-elliptic relative equilibria, negative elliptic relative equilibria, negative hyperbolic relative equilibria, negative elliptic-hyperbolic relative equilibria and higher dimensional versions thereoff) of the $n$-body problem in $\mathbb{S}^{k}$ and $\mathbb{H}^{k}$ and that their proofs do not rely on specific properties of the matrix $T$.
\end{remark}
\begin{remark}
  The restriction that the relative equilibria on the unit sphere lie in $R_{\epsilon, T,m_{1},...,m_{n}}$ is needed to potentially exclude sequences of relative equilibria in $\mathbb{S}^{k}$ that in the limit can show antipodal behaviour, i.e. have sequences of point masses $\{Q_{ip}\}_{p=1}^{\infty}$, $\{Q_{jp}\}_{p=1}^{\infty}$ for which $\lim\limits_{p\rightarrow\infty}\langle Q_{ip},Q_{jp}\rangle=-1$. Proving, or disproving the existence of such sequences could lead to valuable information about the geometry of the $n$-body problem on the unit sphere, but unfortunately lies beyond the scope of this paper.
\end{remark}
We will now prove Theorem~\ref{limits distinct} in section~\ref{Proof of the Main Theorem}, then prove Corollary~\ref{compact} in section~\ref{Proof of the Main Corollary}, after which we will prove Theorem~\ref{Main Theorem Curved Case 1} and Theorem~\ref{Main Theorem Curved Case 2} in section~\ref{Proof of the curved case}.
\section{Proof of Theorem~\ref{limits distinct}}\label{Proof of the Main Theorem}
\begin{proof}
  We will prove this theorem for the case that $\lim\limits_{x\downarrow 0}f(x)=+\infty$ and give an argument at the end of the proof how the argument of the proof can also be used if $\lim\limits_{x\downarrow 0}f(x)=-\infty$.

  Assume that Theorem~\ref{limits distinct} is false. Then there exist sequences $\{Q_{ir}\}_{r=1}^{\infty}$ and relative equilibria $q_{ir}(t)=T(t)Q_{ir}$, $i\in\{1,...,n\}$ for which we may assume, if we renumber the $Q_{ir}$ in terms of $i$ and take subsequences if necessary, the following:
  \begin{itemize}
    \item[1. ] There exist sequences $\{Q_{1r}\}_{r=1}^{\infty}$,...,$\{Q_{lr}\}_{r=1}^{\infty}$, $l\leq n$ such that $\|Q_{ir}-Q_{jr}\|$ goes to zero for $r$ going to infinity if $i$, $j\in\{1,...,l\}$, $2\leq l \geq n$, .
    \item[2. ] $\|Q_{ir}-Q_{jr}\|$ does not go to zero for $r$ going to infinity if $i\in\{1,...,l\}$ and $j\in\{l+1,...,n\}$.
    \item[3. ] $\|Q_{1r}-Q_{2r}\|\geq\|Q_{ir}-Q_{jr}\|$ for all $i$, $j\in\{1,...,l\}$, for all $r\in\mathbb{N}$.
  \end{itemize}
  Write $T(t)^{-1}\ddot{T}(t)=-\mathbf{A}$. Then inserting $q_{ir}(t)=T(t)Q_{i}$, $i\in\{1,...,n\}$ into (\ref{Equations of motion}), using that for any $x\in\mathbb{R}^{k}$ $\|T(t)x\|=\|x\|$ and multiplying both sides of (\ref{Equations of motion}) with $T(t)^{-1}$, gives
  \begin{align*}
    -\mathbf{A}Q_{ir}=\sum\limits_{j=1,\textrm{ }j\neq i}^{n}m_{j}(Q_{jr}-Q_{ir})f(\|Q_{jr}-Q_{ir}\|)
  \end{align*}
  and consequently
  \begin{align}\label{q1}
    -\mathbf{A}Q_{1r}=\sum\limits_{j=2}^{n}m_{j}(Q_{jr}-Q_{1r})f(\|Q_{jr}-Q_{1r}\|)
  \end{align}
  and
  \begin{align}\label{q2}
    -\mathbf{A}Q_{2r}=\sum\limits_{j=1,\textrm{ }j\neq 2}^{n}m_{j}(Q_{jr}-Q_{2r})f(\|Q_{jr}-Q_{2r}\|).
  \end{align}
  Subtracting (\ref{q1}) from (\ref{q2}) gives
  \begin{align}
    \mathbf{A}(Q_{1r}-Q_{2r})&=(m_{1}+m_{2})(Q_{1r}-Q_{2r})f(\|Q_{1r}-Q_{2r}\|)\nonumber\\
    &+\sum\limits_{j=3}^{l}m_{j}\left((Q_{1r}-Q_{jr})f(\|Q_{1r}-Q_{jr}\|)+(Q_{jr}-Q_{2r})f(\|Q_{jr}-Q_{2r}\|)\right)\nonumber\\
    &+\sum\limits_{j=l+1}^{n}m_{j}\left((Q_{1r}-Q_{jr})f(\|Q_{1r}-Q_{jr}\|)+(Q_{jr}-Q_{2r})f(\|Q_{jr}-Q_{2r}\|)\right)\label{go time}
  \end{align}
  Note that for any $j\in\{1,...,l\}$ the vectors $Q_{1r}-Q_{2r}$, $Q_{1r}-Q_{jr}$ and $Q_{jr}-Q_{2r}$ either form a triangle with $\|Q_{1r}-Q_{2r}\|$ the length of its longest side, or the three of them align, meaning the angles between them are zero. Consequently, the angle between $Q_{1r}-Q_{2r}$ and $Q_{1r}-Q_{jr}$ and the angle between $Q_{1r}-Q_{2r}$ and $Q_{jr}-Q_{2r}$ is smaller than $\frac{\pi}{2}$ and thus \begin{align}\label{1>}\lim\limits_{r\rightarrow\infty}\langle Q_{1r}-Q_{jr},Q_{1r}-Q_{2r}\rangle\geq 0\textrm{ and }\lim\limits_{r\rightarrow\infty}\langle Q_{jr}-Q_{2r},Q_{1r}-Q_{2r}\rangle\geq 0.\end{align} Also note that
  \begin{align}
    &\sum\limits_{j=l+1}^{n}m_{j}\left((Q_{1r}-Q_{jr})f(\|Q_{1r}-Q_{jr}\|)+(Q_{jr}-Q_{2r})f(\|Q_{jr}-Q_{2r}\|)\right)\nonumber\\
    &=\sum\limits_{j=l+1}^{n}m_{j}\left((Q_{1r}-Q_{2r}+Q_{2r}-Q_{jr})f(\|Q_{1r}-Q_{jr}\|)+(Q_{jr}-Q_{2r})f(\|Q_{jr}-Q_{2r}\|)\right)\nonumber\\
    &=\sum\limits_{j=l+1}^{n}m_{j}\left((Q_{1r}-Q_{2r})f(\|Q_{1r}-Q_{jr}\|)\right.\nonumber\\
    &\left.+(Q_{jr}-Q_{2r})\left(f(\|Q_{jr}-Q_{2r}\|)-f(\|Q_{1r}-Q_{jr}\|)\right)\right)\label{derivative1}
  \end{align}
  Let
  \begin{align*}
    B_{jr}=\begin{cases}
      \frac{f(\|Q_{jr}-Q_{2r}\|)-f(\|Q_{1r}-Q_{jr}\|)}{\|Q_{jr}-Q_{2r}\|-\|Q_{1r}-Q_{jr}\|} \cdot\frac{\|Q_{jr}-Q_{2r}\|-\|Q_{1r}-Q_{jr}\|}{\|Q_{1r}-Q_{2r}\|}\textrm{ if }\|Q_{jr}-Q_{2r}\|\neq\|Q_{1r}-Q_{jr}\| \\
      0\textrm{ if }\|Q_{jr}-Q_{2r}\|=\|Q_{1r}-Q_{jr}\|
    \end{cases}
  \end{align*}
  Taking inner products on both sides of (\ref{go time}) with $\frac{Q_{1r}-Q_{2r}}{\|Q_{1r}-Q_{2r}\|^{2}}$, using (\ref{1>}), (\ref{derivative1}) and using that for $r$ large enough and $i$, $j\in\{1,...,l\}$ we have that $f(\|Q_{ir}-Q_{jr}\|)\geq 0$, gives for $r$ large enough that
  \begin{align*}
    &\left(\frac{\langle\mathbf{A}(Q_{1r}-Q_{2r}),Q_{1r}-Q_{2r}\rangle}{\|Q_{1r}-Q_{2r}\|^{2}}\right)\geq(m_{1}+m_{2})f(\|Q_{1r}-Q_{2r}\|)+0\\
    &+\sum\limits_{j=l+1}^{n}m_{j}f(\|Q_{1r}-Q_{jr}\|)+\sum\limits_{j=l+1}^{n}m_{j}\langle Q_{jr}-Q_{2r},\frac{Q_{1r}-Q_{2r}}{\|Q_{1r}-Q_{2r}\|}\rangle B_{jr}.
  \end{align*}
  Note that  $\left|\frac{\langle\mathbf{A}(Q_{1r}-Q_{2r}),Q_{1r}-Q_{2r}\rangle}{\|Q_{1r}-Q_{2r}\|^{2}}\right|\leq\frac{\left\|\mathbf{A}(Q_{1r}-Q_{2r})\right\|\left\|Q_{1r}-Q_{2r}\right\|}{\|Q_{1r}-Q_{2r}\|^{2}}\leq c_{2}$, that $\frac{Q_{1r}-Q_{2r}}{\|Q_{1r}-Q_{2r}\|}$ is bounded in norm, that for $\lim\limits_{r\rightarrow\infty}Q_{1r}=\lim\limits_{r\rightarrow\infty}Q_{2r}=Q_{1}$ and $\lim\limits_{r\rightarrow\infty}Q_{jr}=Q_{j}$
  \begin{align*}
    \lim\limits_{r\rightarrow\infty}B_{jr}=\begin{cases}\lim\limits_{r\rightarrow\infty}f'(\|Q_{1r}-Q_{jr}\|)\frac{\|Q_{jr}-Q_{2r}\|-\|Q_{1r}-Q_{jr}\|}{\|Q_{1r}-Q_{2r}\|}\textrm{ if }\|Q_{jr}-Q_{2r}\|\neq\|Q_{1r}-Q_{jr}\|\\ 0\textrm{ if }\|Q_{jr}-Q_{2r}\|=\|Q_{1r}-Q_{jr}\|\end{cases},
  \end{align*}
  which is bounded by construction as by the triangle inequality,
  \begin{align*}
    \left|\frac{\|Q_{jr}-Q_{2r}\|-\|Q_{1r}-Q_{jr}\|}{\|Q_{1r}-Q_{2r}\|}\right|=\frac{|\|Q_{jr}-Q_{2r}\|-\|Q_{1r}-Q_{jr}\||}{\|Q_{1r}-Q_{2r}\|}\leq\frac{\|Q_{1r}-Q_{2r}\|}{\|Q_{1r}-Q_{2r}\|}=1,
  \end{align*}
  which means that we have that for some $B>0$
  \begin{align*}
    \lim\limits_{r\rightarrow\infty}\left|\sum\limits_{j=l+1}^{n}m_{j}f(\|Q_{1r}-Q_{jr}\|)+\sum\limits_{j=l+1}^{n}m_{j}\langle Q_{jr}-Q_{2r},\frac{Q_{1r}-Q_{2r}}{\|Q_{1r}-Q_{2r}\|}\rangle B_{jr}\right|<B
  \end{align*}
  for all $r\in\mathbb{N}$ and thus that
  \begin{align*}
    c_{2}\geq\lim\limits_{r\rightarrow\infty}(m_{1}+m_{2})f(\|Q_{1r}-Q_{2r}\|)-B=\infty,
  \end{align*}
  which is a contradiction. If $\lim\limits_{x\downarrow 0}f(x)=-\infty$, we can define $g=-f$, rewrite everything in terms of $g$ and repeat the proof of our theorem using $g$ instead of $f$. This completes our proof.
\end{proof}
\section{Proof of Corollary~\ref{compact}}~\label{Proof of the Main Corollary}
\begin{proof}
  Assume the contrary to be true. Then there exist sequences $\{Q_{ir}\}_{r=1}^{\infty}$, $i\in\{1,...,n\}$ for which $q_{ir}(t)=T(t)Q_{ir}$ define relative equilibrium solutions of (\ref{Equations of motion}) and for which there has to be at least one sequence $\{Q_{ir}\}_{r=1}^{\infty}$ that is unbounded. Taking subsequences and renumbering the $Q_{ir}$ in terms of $i$ if necessary, we may assume that $\{Q_{1r}\}_{r=1}^{\infty}$ is unbounded. Then by (\ref{q1}),
  \begin{align}\label{Compactness k}
    \mathbf{A}Q_{1r}=\sum\limits_{j=2}^{n}m_{j}(Q_{1r}-Q_{jr})f(\|Q_{1r}-Q_{jr}\|)
  \end{align}
  with $\mathbf{A}=-T(t)^{-1}\ddot{T}(t)$. As the left-hand side of (\ref{Compactness k}) is unbounded, the right-hand side must be unbounded as well, which means that there must be $j\in\{2,...,n\}$ for which \begin{align*}m_{j}(Q_{1r}-Q_{jr})f(\|Q_{1r}-Q_{jr}\|)\end{align*} is unbounded if we let $r$ go to infinity. But as \begin{align*}\left\|m_{j}(Q_{1r}-Q_{jr})f(\|Q_{1r}-Q_{jr}\|)\right\|=m_{j}\|Q_{1r}-Q_{jr}\|f(\|Q_{1r}-Q_{jr}\|),\end{align*} that means that $\|Q_{1r}-Q_{jr}\|$ goes to zero for $r$ going to infinity, which is impossible by Theorem~\ref{limits distinct}. This completes the proof.
\end{proof}
\section{Proof of Theorem~\ref{Main Theorem Curved Case 1} and Theorem~\ref{Main Theorem Curved Case 2}}~\label{Proof of the curved case}
\begin{proof}
  We will prove Theorem~\ref{Main Theorem Curved Case 1} and Theorem~\ref{Main Theorem Curved Case 2} using the same argument:

  Assume that either of the theorems is incorrect. Then there exist sequences of relative equilibria $\{q_{ir}(t)\}_{r=1}^{\infty}=\{T(t)Q_{ir}\}_{r=1}^{\infty}$, $i\in\{1,...,n\}$, for which, renumbering the $q_{ir}$ in terms of $i$ if necessary, there exists an $l\in\{1,...,n\}$ such that
  \begin{itemize}
    \item[1. ] $\|Q_{ir}-Q_{jr}\|$ goes to zero for $r$ going to infinity if $i$, $j\in\{1,...,l\}$.
    \item[2. ] $\|Q_{ir}-Q_{jr}\|$ does not go to zero for $r$ going to infinity if $i\in\{1,...,l\}$ and $j\in\{l+1,...,n\}$.
  \end{itemize}
  Inserting these solutions into (\ref{EquationsOfMotion Curved}), using that for $x$,$y\in\mathbf{M}_{\sigma}^{k}$ $(T(t)x)\odot (T(t)y)=x\odot y$, subsequently multiplying both sides of the equation from the left with $T(t)^{-1}$ and defining $T(t)^{-1}\ddot{T}(t)Q_{ir}+\sigma(\dot{q}_{ir}\odot\dot{q}_{ir})Q_{ir}=B_{\sigma ir}$ gives
  \begin{align}\label{Quick}
    B_{\sigma ir}=\sum\limits_{j=1,\textrm{ }j\neq i}^{n}\frac{m_{j}(Q_{jr}-\sigma(Q_{ir}\odot Q_{jr})Q_{ir})}{(\sigma-\sigma(Q_{ir}\odot Q_{jr})^{2})^{\frac{3}{2}}}.
  \end{align}
  Multiplying both sides of (\ref{Quick}) with $m_{i}$ and summing over $i$ from $1$ to $l$ then gives
  \begin{align}\label{Quick2}
    \sum\limits_{i=1}^{l}m_{i}B_{\sigma ir}=\sum\limits_{i=1}^{l}\sum\limits_{j=1,\textrm{ }j\neq i}^{n}\frac{m_{i}m_{j}(Q_{jr}-\sigma(Q_{ir}\odot Q_{jr})Q_{ir})}{(\sigma-\sigma(Q_{ir}\odot Q_{jr})^{2})^{\frac{3}{2}}}.
  \end{align}
  Note that
  \begin{align*}
    &2\sum\limits_{i=1}^{l}\sum\limits_{j=1,\textrm{ }j\neq i}^{l}\frac{m_{i}m_{j}(Q_{jr}-\sigma(Q_{ir}\odot Q_{jr})Q_{ir})}{(\sigma-\sigma(Q_{ir}\odot Q_{jr})^{2})^{\frac{3}{2}}}\\
    &=\left(\sum\limits_{i=1}^{l}\sum\limits_{j=1,\textrm{ }j\neq i}^{l}\frac{m_{i}m_{j}(Q_{jr}-\sigma(Q_{ir}\odot Q_{jr})Q_{ir})}{(\sigma-\sigma(Q_{ir}\odot Q_{jr})^{2})^{\frac{3}{2}}}+\sum\limits_{j=1}^{l}\sum\limits_{i=1,\textrm{ }i\neq j}^{l}\frac{m_{j}m_{i}(Q_{ir}-\sigma(Q_{jr}\odot Q_{ir})Q_{jr})}{(\sigma-\sigma(Q_{jr}\odot Q_{ir})^{2})^{\frac{3}{2}}}\right),
  \end{align*}
  which, collecting all terms that are multiplied with $m_{i}m_{j}$, can be rewritten as
  \begin{align*}
    &2\sum\limits_{i=1}^{l}\sum\limits_{j=1,\textrm{ }j\neq i}^{l}\frac{m_{i}m_{j}(Q_{jr}-\sigma(Q_{ir}\odot Q_{jr})Q_{ir})}{(\sigma-\sigma(Q_{ir}\odot Q_{jr})^{2})^{\frac{3}{2}}}\nonumber\\
    &=\sum\limits_{i=1}^{l}\sum\limits_{j=1,\textrm{ }j\neq i}^{l}\frac{m_{i}m_{j}(Q_{ir}+Q_{jr})(1-\sigma(Q_{ir}\odot Q_{jr}))}{(\sigma-\sigma(Q_{ir}\odot Q_{jr})^{2})^{\frac{3}{2}}}\\
    &=\sum\limits_{i=1}^{l}\sum\limits_{j=1,\textrm{ }j\neq i}^{l}\frac{m_{i}m_{j}(Q_{ir}+Q_{jr})(1-\sigma(Q_{ir}\odot Q_{jr}))}{(\sigma-\sigma(Q_{ir}\odot Q_{jr})^{2})^{\frac{3}{2}}}\cdot\frac{\sigma(1+\sigma(Q_{ir}\odot Q_{jr}))}{\sigma(1+\sigma(Q_{ir}\odot Q_{jr}))},
  \end{align*}
  giving
  \begin{align}
    &2\sum\limits_{i=1}^{l}\sum\limits_{j=1,\textrm{ }j\neq i}^{l}\frac{m_{i}m_{j}(Q_{jr}-\sigma(Q_{ir}\odot Q_{jr})Q_{ir})}{(\sigma-\sigma(Q_{ir}\odot Q_{jr})^{2})^{\frac{3}{2}}}\nonumber\\
    &=\sum\limits_{i=1}^{l}\sum\limits_{j=1,\textrm{ }j\neq i}^{l}\frac{m_{i}m_{j}(Q_{ir}+Q_{jr})(\sigma-\sigma(Q_{ir}\odot Q_{jr})^{2})}{(\sigma-\sigma(Q_{ir}\odot Q_{jr})^{2})^{\frac{3}{2}}\sigma(1+Q_{ir}\odot Q_{jr})}\nonumber\\
    &=\sum\limits_{i=1}^{l}\sum\limits_{j=1,\textrm{ }j\neq i}^{l}\frac{m_{i}m_{j}(Q_{ir}+Q_{jr})}{(\sigma-\sigma(Q_{ir}\odot Q_{jr})^{2})^{\frac{1}{2}}\sigma(1+Q_{ir}\odot Q_{jr})}. \label{double sum trick curved}
  \end{align}
  Let
  \begin{align}\label{B2sigma ir}
    B_{2\sigma r}=\sum\limits_{i=1}^{l}\sum\limits_{j=l+1}^{n}\frac{m_{i}m_{j}(Q_{jr}-\sigma(Q_{ir}\odot Q_{jr})Q_{ir})}{(\sigma-\sigma(Q_{ir}\odot Q_{jr})^{2})^{\frac{3}{2}}}.
  \end{align}
  Inserting (\ref{double sum trick curved}) and (\ref{B2sigma ir}) into (\ref{Quick2}) gives
  \begin{align*}
    \sum\limits_{i=1}^{l}m_{i}B_{\sigma ir}Q_{ir}-B_{2\sigma r}=\frac{1}{2}\sigma\sum\limits_{i=1}^{l}\sum\limits_{j=1,\textrm{ }j\neq i}^{l}\frac{m_{i}m_{j}(Q_{ir}+Q_{jr})}{(\sigma-\sigma(Q_{ir}\odot Q_{jr})^{2})^{\frac{1}{2}}(1+Q_{ir}\odot Q_{jr})},
  \end{align*}
  which means that
  \begin{align}\label{and presto curved}
    \left(\sum\limits_{i=1}^{l}m_{i}B_{\sigma ir}Q_{ir}-B_{2\sigma r}\right)\odot Q_{1r}=\frac{1}{2}\sum\limits_{i=1}^{l}\sum\limits_{j=1,\textrm{ }j\neq i}^{l}\frac{m_{i}m_{j}(Q_{ir}+Q_{jr})\odot Q_{1r}}{(\sigma-\sigma(Q_{ir}\odot Q_{jr})^{2})^{\frac{1}{2}}\sigma(1+Q_{ir}\odot Q_{jr})}.
  \end{align}
  By construction, $\left(\sum\limits_{i=1}^{l}m_{i}B_{\sigma ir}Q_{ir}-B_{2\sigma r}\right)\odot Q_{1r}$ and $\sigma(1+Q_{ir}\odot Q_{jr})>0$ are bounded for $r$ going to infinity and $\lim\limits_{r\rightarrow\infty}Q_{ir}\odot Q_{jr}=1$ for $i$, $j\in\{1,...,l\}$, so by (\ref{and presto curved})
  \begin{align}
    &\lim\limits_{r\rightarrow\infty}\left(\sum\limits_{i=1}^{l}m_{i}B_{\sigma ir}Q_{ir}-B_{2\sigma r}\right)\odot Q_{1r}\nonumber\\
    &=\frac{1}{2}\sum\limits_{i=1}^{l}\sum\limits_{j=1,\textrm{ }j\neq i}^{l}\frac{m_{i}m_{j}(Q_{ir}+Q_{jr})\odot Q_{1r}}{(\sigma-\sigma(Q_{ir}\odot Q_{jr})^{2})^{\frac{1}{2}}\sigma(1+Q_{ir}\odot Q_{jr})}=\infty,\label{done curved}
  \end{align}
  while the left-hand side of (\ref{done curved}) is bounded, which is a contradiction. This completes the proof.
\end{proof}

\end{document}